\documentclass[12pt]{article} 
\usepackage{amsmath}
\usepackage{amssymb} 
\usepackage{amscd}

\usepackage{theorem} 
\usepackage[T1]{fontenc}
\usepackage{xspace}

\def\dim{{\mbox{dim}}}

\def\Hom {{\mbox{Hom}}}

\def\cala{{\cal A}}

\def\calc{{\cal C}}

 \def\fracg{{\mathfrak g}} 
   
\def\bbbone{\mbox{\rm 1\hspace {-.6em} l}}

\def\gg{{\mathbf g}}

\newtheorem{theorem}{THEOREM}
 
\newtheorem{proposition}{PROPOSITION}
\newtheorem{corollary}{COROLLARY}

\begin{document}

\baselineskip=0.5cm

\begin{center} 
 \thispagestyle{empty}
{\large\bf YANG-MILLS ALGEBRA}
\end{center} 
\vspace{0.5cm}

\begin{center} Alain CONNES \footnote{Coll\`ege de France, 3 rue
d'Ulm, 75  005 Paris, and \\ I.H.E.S., 35 route de Chartres, 91440 
Bures-sur-Yvette\\
connes$@$ihes.fr}
and 
 Michel DUBOIS-VIOLETTE
\footnote{Laboratoire de Physique Th\'eorique, UMR 8627, Universit\'e Paris XI,
B\^atiment 210, F-91 405 Orsay Cedex, France\\
Michel.Dubois-Violette$@$th.u-psud.fr\\
}
\end{center} \vspace{0,25cm}

\begin{center} \today \end{center}

\vspace {0,25cm}

\begin{abstract}
 Some unexpected properties of the cubic algebra generated by the covariant derivatives of a generic Yang-Mills connection over the $(s+1)$-dimensional pseudo Euclidean space are pointed out. This algebra is Gorenstein and Koszul of global dimension 3 but except for $s=1$ (i.e. in the 2-dimensional case) where it is the universal enveloping algebra
of the Heisenberg Lie algebra and is a cubic Artin-Schelter regular algebra, it fails to be regular in that it has exponential growth.
We give an explicit formula for the
Poincar\'e series of this algebra $\cala$ and for
the dimension in degree $n$ of the graded Lie algebra
of which $\cala$ is the universal enveloping algebra.
 In the 4-dimensional (i.e. $s=3$) Euclidean case, a quotient of this algebra is the quadratic algebra generated by the covariant derivatives of a generic (anti) self-dual connection. This latter algebra is Koszul of global dimension 2 but is not Gorenstein and has exponential growth. It is the universal enveloping algebra of the
graded Lie-algebra which is the semi-direct product of the free
Lie algebra with three generators of degree one by a derivation of
degree one.
\end{abstract}

\vspace{0,25cm}

  \vspace{1cm}

\noindent LPT-ORSAY 02-50

\newpage

\section{Introduction}
 Our aim here is to investigate the properties of a cubic algebra which is behind the Yang-Mills equations in ($s+1$)-dimensional pseudo Euclidean space.\\

Throughout the following, $s$ denotes an integer with $s\geq 1$, the $x^\lambda$ denote the canonical coordinates of $\mathbb R^{s+1}$, the corresponding partial derivatives are denoted by $\partial_\lambda=\frac{\partial}{\partial x^\lambda}$ and the pseudo Euclidean structure on $\mathbb R^{s+1}$ is defined by the scalar products $g_{\mu\nu}$ of the corresponding basis elements. Thus the $g_{\mu\nu}$ for $\mu,\nu\in \{0,1,\dots,s\}$ are the components of an invertible real symmetric matrix and we shall denote by $g^{\mu\nu}$ the matrix elements of the inverse matrix, that is one has $g_{\mu\lambda}g^{\lambda \nu}=\delta^\nu_\mu$.\\

The Yang-Mills equations here are equations for connections on bundles over $\mathbb R^{s+1}$. Any complex vector bundle of rank $p$ over $\mathbb R^{s+1}$ is isomorphic to the trivial bundle $\mathbb C^p\otimes \mathbb R^{s+1}\rightarrow \mathbb R^{s+1}$  and a connection on such a bundle is described by a $M_p(\mathbb C)$-valued one-form $A_\mu dx^\mu$. The corresponding covariant derivatives are given by $\nabla_\lambda=\partial_\lambda+A_\lambda$. The Yang-Mills equations for the $A_\mu$ are compactly expressed in terms of the $\nabla_\mu$ by 
\[
g^{\lambda\mu}[\nabla_\lambda,[\nabla_\mu,\nabla_\nu]]=0
\]
for $\nu\in \{0,\dots,s\}$.\\

It follows that any solution of the Yang-Mills equations carries a representation of the abstract algebra $\cala$ generated by elements $\nabla_\lambda$ with the cubic relations above. Recently there has been a renewal of interest for this algebra in connection with string theory \cite{NN}.\\

The algebra $\cala$ is a homogeneous algebra of degree 3, in short a cubic algebra, so the concepts and technics developed in \cite{BD-VW} and \cite{RB3} are well-suited to analyse its properties. This is the aim of Sections 2 and 3.\\

In the 4-dimensional case (i.e. $s=3$) with the Euclidean structure $g_{\mu\nu}=\delta_{\mu\nu}$ ($\mu,\nu=0,1,2,3$) and the orientation of $\mathbb R^4$ corresponding to the order $0,1,2,3$, a connection can be self-dual that is such that 
\[
[\nabla_0,\nabla_k]=[\nabla_\ell,\nabla_m]
\]
or anti-self-dual that is such that
\[
[\nabla_0,\nabla_k]=-[\nabla_\ell,\nabla_m]
\]
for any cyclic permutation $(k, \ell, m)$ of ($1,2,3$). A self-dual or anti-self-dual connection satisfies the Yang-Mills equations which means that the cubic algebra $\cala$ admits as quotients the quadratic algebras $\cala^{(+)}$ and $\cala^{(-)}$ generated by elements $\nabla_\lambda$ ($\lambda=0,1,2,3$) satisfying the corresponding quadratic relations above. These algebras are investigated in Section 4.\\

The algebras investigated here belong to the class of homogeneous algebras, a class of algebras which play a key role in the noncommutative version of algebraic geometry (see e.g. \cite{AS}, \cite{ATVB}, \cite{Kempf}, \cite{SmSt}, \cite{Tvdb}, \cite{RB3}, \cite{AC.MDV}). In this context the canonical generalization of the Koszul complexes of quadratic algebras \cite{Pri}, \cite{YuM2} are $N$-complexes \cite{BD-VW} and the right generalization of the notion of Koszul algebra introduced in \cite{RB3} is closely tied to these $N$-complexes \cite{BD-VW}. It is worth noticing here that the theory of $N$-complexes has independently benefited of recent developments (see e.g. \cite{D-V2}, \cite{D-VH2}, \cite{D-VK}, \cite{D-VT2}, \cite{Kap}, \cite{KW}, \cite{Wam3} and references quoted in these papers).\\

For the convenience of the reader we have added an appendix summarizing a simplified version of some notions used in the paper for homogeneous algebras \cite{BD-VW} such as the Koszul and the Gorenstein properties.

 \section{The cubic Yang-Mills algebra and its dual}
 
 Let $\cala$ be the unital associative $\mathbb C$-algebra generated by the elements $\nabla_\lambda$, $\lambda\in \{0,1,\dots,s\}$ with relations 
\begin{equation}
g^{\lambda\mu}[\nabla_\lambda,[\nabla_\mu,\nabla_\nu]]=0,\>\>\> \nu\in\{0,1,\dots,s\}
\label{eq1}
\end{equation}
By giving the degree 1 to the $\nabla_\lambda$, $\cala$ is a graded algebra which is connected (i.e. $\cala_0=\mathbb C\bbbone$) and generated in degree 1. In fact $\cala$ is the 3-homogeneous  (or cubic) algebra \cite{RB3}, \cite{BD-VW} $\cala=A(E,R)$ with $E=\oplus_\lambda \mathbb C \nabla_\lambda$ and where $R\subset E^{\otimes^3}$ is spanned by the 
\begin{equation}
g^{\lambda\mu}(\nabla_\lambda \otimes \nabla_\mu \otimes \nabla_\nu + \nabla_\nu \otimes \nabla_\lambda \otimes \nabla_\mu-2\nabla_\lambda \otimes \nabla_\nu \otimes \nabla_\mu)
\label{eq2}
\end{equation}
for $\nu \in \{ 0,1,\dots,s\}$. The cubic algebra $\cala$ will be refered to as {\sl the cubic Yang-Mills algebra}. The algebra $\cala$ can be equipped with a structure of $\ast$-algebra if one defines its involution $a\mapsto a^\ast$ by setting
\begin{equation}
\nabla^\ast_\lambda = -\nabla_\lambda
\label{eq3}
\end{equation}
for $\lambda \in \{0,1,\dots,s\}$.\\

The dual $\cala^!$ of $\cala$ is by definition \cite{BD-VW} the 3-homogeneous algebra $\cala^!=A(E^\ast,R^\perp)$ where $E^\ast=\oplus_\lambda \mathbb C \theta^\lambda$ is the dual vector space of $E$ with $\theta^\lambda(\nabla_\mu)=\delta^\lambda_\mu$ and where $R^\perp\subset E^{\ast \otimes^3}$ is the annihilator of $R$. It is easy to verify that $\cala^!$ is the unital associative $\mathbb C$-algebra generated by the $\theta^\lambda$, $\lambda\in \{0,1,\dots, s\}$ with relations
\begin{equation}
\theta^\lambda\theta^\mu\theta^\nu=\frac{1}{s}(g^{\lambda\mu}\theta^\nu+g^{\mu\nu}\theta^\lambda -2g^{\lambda\nu}\theta^\mu)\gg 
\label{eq4}
\end{equation}
where $\gg=g_{\alpha\beta}\theta^\alpha\theta^\beta$. By contraction of Relation (\ref{eq4}) with $g_{\lambda\mu}$ one obtains
\begin{equation}
\gg \theta^\nu=\theta^\nu \gg
\label{eq5}
\end{equation}
which means that $\gg\>  (\in \cala^!_2)$ is central. Using the bilinear scalar product on $E^\ast$ defined by $(a\vert b)=g^{\lambda\mu}a_\lambda b_\mu$ for $a=a_\lambda\theta^\lambda$, $b=b_\lambda\theta^\lambda\in E^\ast=\cala^!_1$, one can rewrite (\ref{eq4}) in the form
\begin{equation}
abc = \frac{1}{s}((a\vert b)c+(b\vert c)a-2(a\vert c)b)\gg
\label{eq6}
\end{equation}
for any $a,b,c\in E^\ast = \cala^!_1$. It follows in particular that one has $abc=cba$ and $a^3=0$ for any $a,b,c\in\cala^!_1$.

\begin{proposition}
One has $\cala^!_0=\mathbb C\bbbone$, $\cala^!_1=\oplus_\lambda\mathbb C\theta^\lambda$, $\cala^!_2=\oplus_{\mu\nu} \mathbb C \theta^\mu\theta^\nu$, $\cala^!_3=\oplus_\lambda \mathbb C \theta^\lambda \gg$, $\cala^!_4=\mathbb C \gg^2$ and $\cala^!_n=0$ for $n\geq 5$.
\end{proposition}

The first 3 equalities follow from the definition, the last 3 ones are easily established for instance by using (\ref{eq4}) in a basis ($\theta^\lambda$) where $g^{\mu\nu}$ is diagonal. Thus, one has $\dim(\cala^!_0)=\dim(\cala^!_4)=1$, $\dim(\cala^!_1)=\dim(\cala^!_3)=s+1$, $\dim(\cala^!_2)=(s+1)^2$ and $\dim(\cala^!_n)=0$ otherwise.

\section{Homological properties}

As explained in \cite{BD-VW} to the cubic algebra $\cala$ is associated the chain 3-complex of left $\cala$-modules $K(\cala)$ defined by
\[
0\rightarrow \cala\otimes (\cala^!_4)^\ast\stackrel{d}{\rightarrow}\cala \otimes (\cala^!_3)^\ast \stackrel{d}{\rightarrow} \cala\otimes (\cala^!_2)^\ast \stackrel{d}{\rightarrow} \cala\otimes (\cala^!_1)^\ast \stackrel{d}{\rightarrow} \cala\rightarrow 0
\]
with the $(\cala^!_n)^\ast$ given here by $(\cala^!_1)^\ast=E$, $(\cala^!_2)^\ast=E^{\otimes^2}$, $(\cala^!_3)^\ast =R\subset E^{\otimes^3}$, $(\cala^!_4)^\ast=(E\otimes R)\cap (R\otimes E)\subset E^{\otimes^4}$ and where $d$ is induced by 
\[
a\otimes(e_1\otimes\dots \otimes e_n)\mapsto ae_1\otimes (e_2\otimes\dots \otimes e_n)
\]
for $a\in \cala$ and $e_i\in \cala_1=E$. From $K(\cala)$ one extracts the chain complex of left $\cala$-module $C_{2,0}$ defined by
\[
0\rightarrow \cala\otimes (\cala^!_4)^\ast \stackrel{d}{\rightarrow} \cala\otimes (\cala^!_3)^\ast \stackrel{d^2}{\rightarrow} \cala\otimes (\cala^!_1)^\ast \stackrel{d}{\rightarrow} \cala\rightarrow 0
\]
which, as pointed out in \cite{BD-VW}, coincides with the Koszul complex of \cite{RB3}. We shall show that this complex is acyclic in positive degrees $n\geq 1$ which means that $\cala$ is Koszul in the sense of \cite{RB3}.

\begin{theorem}
The cubic Yang-Mills algebra $\cala$ is Koszul of global dimension 3 and is Gorenstein.
\end{theorem}

\noindent \underbar{Proof}. In view of the results of last section on the $\cala^!_n$, one has the isomorphisms of left $\cala$-modules $\cala\otimes (\cala^!_4)^\ast \simeq \cala$ and $\cala\otimes (\cala^!_3)^\ast\simeq \cala\otimes (\cala^!_1)^\ast \simeq \cala^{s+1}=\underbrace{(\cala,\dots,\cala)}_{s+1}$ and with the corresponding identifications the complex $C_{2,0}$ reads
\begin{equation}
0\rightarrow \cala\stackrel{\nabla^t}{\rightarrow}\cala^{s+1} \stackrel{M}{\rightarrow} \cala^{s+1} \stackrel{\nabla}{\rightarrow} \cala\rightarrow 0
\label{eq7}
\end{equation}
where $\nabla$ is the column matrix with components $\nabla_\lambda$ while $\nabla^t$ is its transposed that is $\nabla^t=(\nabla_0,\dots,\nabla_s)$ and where $M$ is the $(s+1)\times (s+1)$-matrix with components
\begin{equation}
M^{\mu\nu}=(g^{\mu\nu}g^{\alpha\beta}+g^{\mu\alpha}g^{\nu\beta}-2g^{\mu\beta}g^{\nu\alpha})\nabla_\alpha\nabla_\beta
\label{eq8}
\end{equation}
The relations (\ref{eq1}) defining $\cala$ read $M\nabla=0$ or equivalently $\nabla^tM=0$ as easily verified. This implies that one has the exact sequence of left $\cala$-modules
\begin{equation}
0\rightarrow \cala\stackrel{\nabla^t}{\rightarrow}\cala^{s+1}\stackrel{M}{\rightarrow} \cala^{s+1} \stackrel{\nabla}{\rightarrow} \cala \stackrel{\varepsilon}{\rightarrow} \mathbb C\rightarrow 0
\label{eq9}
\end{equation}
where $\varepsilon$ is the character of $\cala$ defined by the projection on degree 0; (exactness of $\cala^{s+1}\stackrel{M}{\rightarrow} \cala^{s+1}\stackrel{\nabla}{\rightarrow} \cala \rightarrow \mathbb C \rightarrow 0$ follows from the definitions). This implies of course the Koszulity of $\cala$ but also the Gorenstein property by using the obvious symmetry by transposition. The sequence (\ref{eq9}) is in fact a minimal projective resolution of the trivial $\cala$-module $\mathbb C$ in the graded category of $\cala$ which implies by standard arguments that $\cala$ has global dimension 3, \cite{AS}, (see e.g. in \cite{RB4}).$\square$\\

Let us remind that the Gorenstein property means here (see the appendix) that the cochain complex of right $\cala$-modules obtained by applying the functor $\Hom_{\cala}(\bullet, \cala)$ (i.e. by duality) to each left $\cala$-module of the Koszul complex (\ref{eq7}) gives a resolution of the trivial right $\cala$-module $\mathbb C$. In other words this means that this cochain complex $C'_{2,0}$ dual of the Koszul complex $C_{2,0}$ is acyclic in degrees smaller than the dimension which is here 3 (i.e. $H^k(C'_{2,0})=0$, $\forall k<3$) and such that its cohomology in degree 3 coincides with $\mathbb C$ (i.e. $H^3(C'_{2,0})\simeq \mathbb C$),  that is $H(C'_{2,0})=H^3(C'_{2,0})\simeq \mathbb C$.\\

By using the results of \cite{RB3}, (Section 5 and the erratum), one obtains the following corollary.

\begin{corollary}
The left $\cala\otimes \cala^{opp}$-module $\cala$ admits a minimal projective resolution of the form
\[
0\rightarrow \cala\otimes\cala^{opp}\rightarrow \cala^{s+1}\otimes \cala^{opp}\rightarrow \cala^{s+1}\otimes \cala^{opp}\rightarrow \cala\otimes \cala^{opp}\rightarrow \cala\rightarrow 0
\]
in the graded category of $\cala\otimes \cala^{opp}$.
\end{corollary}
This corollary opens the way for the computation of the Hochschild and of the cyclic (co)homology of $\cala$ \cite{connes:02}, \cite{connes:03} and implies in particular the following result.
\begin{corollary}
The cubic Yang-Mills algebra $\cala$ has Hochschild dimension 3.
\end{corollary}
Using the Koszul property and the result of Section 2 on the dimension of the $\cala^!_n$ one obtains the following description of the Poincar\'e series 
$P_\cala(t)=\sum_{n\geq 0} \dim (\cala_n) t^n$ of $\cala$
\begin{corollary}
One has $P_\cala(t)=\frac{1}{(1-t^2)(1-(s+1)t+t^2)}$
\end{corollary}

\noindent \underbar{Proof}. Since $\cala$ is a cubic Koszul algebra one has by applying Corollary~1 of \cite{D-VP} with $N=3$
\[
P_\cala (t) Q_\cala(t) =1
\]
where
\[
Q_\cala(t) = \sum_n (\dim(\cala^!_{3n})t^{3n} -\dim(\cala^!_{3n+1})t^{3n+1})
\]
that is here where
\[
Q_\cala(t)=1-(s+1)t+(s+1)t^3-t^4=(1-t^2)(1-(s+1)t+t^2)
\]
in view of Proposition 1. $\square$\\

From this result, one sees that $\cala$ has exponential growth for $s\geq 2$ while it is of polynomial growth for $s=1$ but in this latter case it is a well-known cubic Artin-Schelter regular algebra \cite{AS}, \cite{ATVB} which is the universal enveloping algebra of the Heisenberg Lie algebra.\\

By construction the Yang-Mills algebra $\cala$ is the universal enveloping algebra of the graded Lie algebra $\fracg=\oplus_{j\geq 1}\fracg_j$ with $s+1$ generators in degree 1 and the same defining relations (\ref{eq1}). The Poincar\'e-Birkhoff-Witt theorem thus gives the following relation between the dimensions $N_j$ of the $\fracg_j$ and the Poincar\'e series $P_\cala(t)$
\[
\prod_j (\frac{1}{1-t^j})^{N_j}=P_\cala(t).
\]
Therefore one gets the explicit formula for $j>2$
\[
N_j=\frac{1}{j} \sum_k \mu(\frac{j}{k}) (t^k_1+t^k_2)
\]
where $t_1$ and $t_2$ are the roots of the equation $t^2-(s+1) t+1=0$ and $\mu(x)$ is the M\"obius function which vanishes unless $x$ is an integer of the form $p_1\dots p_r$ where the $p_i$ are distinct primes  in which case its value is $(-1)^r$. In the interesting case $s+1=4$, this gives the following explicit first
values
for the dimension 
$N_j$ with $j=1,2,....$
\[
\begin{array}{l}
\{4, 6,  16, 45, 144, 440, 1440, 4680, 15600, 52344, 177840, 608160, 2095920,\\
7262640, 25300032, 88517520,310927680, 1095923400, 3874804560,\\ 13737892896,
48829153920, 173949483240, 620963048160, 2220904271040, \\7956987570576,
28553731537320, 102617166646800, 369294887482560,\\
 1330702217420400,
4800706662984672,....\}.
\end{array}
\]

\section{The quadratic (anti) self-duality algebra}

In this section, $s=3$ and $g_{\mu\nu}=\delta_{\mu\nu}$ is the canonical Euclidean metric of $\mathbb R^4$. In this case the cubic algebra $\cala$ admits two non trivial quotients $\cala^{(+)}$ and   
 $\cala^{(-)}$  which are quadratic algebras. Let $\varepsilon=\pm$ and let $\cala^{(\varepsilon)}$ be the unital associative $\mathbb C$-algebra generated by the elements $\nabla_\lambda$, $\lambda\in\{0,1,2,3\}$ with relations
\[
[\nabla_0,\nabla_k]=\varepsilon[\nabla_\ell,\nabla_m]
\]
for any cyclic permutation ($k,\ell,m$) of $(1,2,3)$. The quadratic algebra $\cala^{(+)}$ will be refered to as {\sl the quadratic self-duality algebra} whereas $\cala^{(-)}$ will be refered to as {\sl the quadratic anti-self-duality algebra}. One exchanges $\cala^{(+)}$ and $\cala^{(-)}$ by changing the orientation of $\mathbb R^4$. In the following we shall only consider $\cala^{(+)}$, that is the algebra generated by the elements $\nabla_0,\nabla_1,\nabla_2,\nabla_3$ with relations 
\begin{equation}
[\nabla_0,\nabla_k]=[\nabla_\ell,\nabla_m]
\end{equation}
for any cyclic permutation ($k,\ell,m$) of $(1,2,3)$.\\

The dual quadratic algebra $\cala^{(+)!}$ of $\cala^{(+)}$ is generated by the elements $\theta^\lambda$, $\lambda\in\{0,1,2,3\}$ (the dual basis of the $\nabla_\lambda$ as basis of $\mathbb C^4$) with the relations $\theta^\mu \theta^\nu + \theta^\nu \theta^\mu=0,\> \> \> \forall \mu,\nu\in \{0,1,2,3\}$ and $[\theta^0,\theta^k]+[\theta^\ell,\theta^m]=0$ for any cyclic permutation $(k,\ell, m)$ of $(1,2,3)$, that is with the relations 
\begin{equation}
\theta^\lambda \theta^\mu + \frac{1}{2} \sum_{\nu,\rho} \epsilon^{\lambda \muÊ\nu\rho} \theta^\nu \theta^\rho
\label{eq11}
\end{equation}
where $\epsilon^{\lambda \muÊ\nu\rho}$ is completely antisymmetric with $\epsilon^{0123}=1$.\\
Relations (\ref{eq11}) implies that $\theta^\alpha\theta^\beta\theta^\gamma=0$ for any $\alpha, \beta,\gamma\in \{ 0,1,2,3\}$ so one has the following result.

\begin{proposition}
The homogeneous components of $\cala^{(+)!}$ are given by
\[
\cala^{(+)!}_0=\mathbb C\bbbone,\> \>  \cala^{(+)!}_1=\oplus_{\lambda=0}^3\mathbb C\theta^\lambda,\>\> 
\cala^{(+)!}_2=\oplus_{k=1}^3 \mathbb C \theta^0\theta^k
\]
and $\cala^{(+)!}_n=0$ for $n\geq 3$.
\end{proposition}

The algebra $\cala^{(+)}$ is a quadratic algebra so one has the usual Koszul complexe $K(\cala^{(+)})$, \cite{YuM2}, which reduces here to
\[
0\rightarrow \cala^{(+)}\otimes \left( \cala^{(+)!}_2\right)^\ast \stackrel{d}{\rightarrow} \cala^{(+)}\otimes \left( \cala^{(+)!}_1\right)^\ast \stackrel{d}{\rightarrow} \cala^{(+)}\rightarrow 0
\]
and by using Proposition 2 to
\begin{equation}
0\rightarrow (\cala^{(+)})^3 \stackrel{N}{\rightarrow}(\cala^{(+)})^4\stackrel{\nabla}{\rightarrow} \cala^{(+)} \rightarrow 0
\label{eq12}
\end{equation}
where $N$ is the $3\times 4$-matrix given by
\begin{equation}
N = \left (
\begin{array}{cccc}
-\nabla_1 & \nabla_0 & \nabla_3 & -\nabla_2\\
-\nabla_2 & -\nabla_3 & \nabla_0 & \nabla_1\\
-\nabla_3 & \nabla_2 & -\nabla_1 & \nabla_0
\end{array}
\right)
\label{eq13}
\end{equation}
and where $\nabla$ is the column matrix with component $\nabla_\lambda$.

\begin{theorem}
The quadratic self-duality algebra $\cala^{(+)}$ is a Koszul algebra of global dimension 2.\
\end{theorem}

\noindent \underbar{Proof}. The sequence of left $\cala^{(+)}$-modules
\begin{equation}
0\rightarrow (\cala^{(+)})^3\stackrel{N}{\rightarrow} (\cala^{(+)})^4\stackrel{\nabla}{\rightarrow} \cala^{(+)} \stackrel{\varepsilon}{\rightarrow} \mathbb C \rightarrow 0
\label{eq14}
\end{equation}
is exact. In fact exactness of $(\cala^{(+)})^3\stackrel{N}{\rightarrow} \dots \stackrel{\varepsilon}{\rightarrow} \mathbb C\rightarrow 0$
follows from the definition and the injectivity of $(\cala^{(+)})^3\stackrel{N}{\rightarrow} (\cala^{(+)})^4$ is easy to verify.$\square$

\begin{corollary}
The left $\cala^{(+)}\otimes \cala^{(+)opp}$-module $\cala^{(+)}$ admits a minimal projective resolution of the form
\[
0\rightarrow (\cala^{(+)})^3\otimes \cala^{(+)opp}\rightarrow 
(\cala^{(+)})^4\otimes \cala^{(+)opp}\rightarrow
\cala^{(+)}\otimes \cala^{(+)opp}\rightarrow \cala^{(+)}\rightarrow 0
\]
in the graded category of $\cala^{(+)}\otimes \cala^{(+) opp}$.
\end{corollary}
This corollary then implies the following one.

\begin{corollary}
The quadratic self-duality algebra $\cala^{(+)}$ has Hochschild dimension 2.
\end{corollary}
Finally by the Koszul property, one can compute the Poincar\'e series $P_{\cala^{(+)}}$ of $\cala^{(+)}$ which is given by the following corollary.

\begin{corollary}
One has $P_{\cala^{(+)}}(t)=\frac{1}{(1-t)(1-3t)}$
\end{corollary}

From this result one sees that $\cala^{(+)}$ has exponential growth. On the other hand, the asymmetry of the Koszul complex shows that $\cala^{(+)}$ cannot be Gorenstein.\\

Corollary 6 has the following interpretation. The algebra $\cala^{(+)}$ is the universal enveloping algebra of the semi-direct product of the free Lie algebra $L(\nabla_1,\nabla_2,\nabla_3)$ over 3 generators $\nabla_1,\nabla_2,\nabla_3$ by the derivation $\delta$ given by $\delta(\nabla_k)=[\nabla_\ell,\nabla_m]$ for any cyclic permutation ($k,\ell, m$) of ($1,2,3$). This implies that $\cala^{(+)}$ is itself the cross-product of the tensor algebra $T(\mathbb C^3)=\mathbb C\langle \nabla_1,\nabla_2,\nabla_3\rangle$ with the corresponding derivation (remembering that the algebra $\mathbb C\langle \nabla_1,\nabla_2,\nabla_3\rangle$ is the universal enveloping algebra of the free Lie algebra $L(\nabla_1,\nabla_2,\nabla_3)$ \cite{Bour}). Corollaries 6 follows directly and by remembering that the Hochschild dimension of a tensor algebra is 1 this also implies Corollary 5. It is worth noticing here that this presentation breaks the symmetry $SO(4)$ to retain only the $SO(3)$ of $(\nabla_1,\nabla_2,\nabla_3)$.

 \section*{Acknowledgements} 
 
 We thank Roland Berger for his kind advices and Nikita Nekrasov for renewing our interest on the Yang-Mills algebra.
  
  \section*{Appendix}
  
  Our aim in this appendix is to review in the case of homogeneous algebras some basic notions used in this paper. Many of these notions like the Gorenstein property exist in the more general setting of connected $\mathbb N$-graded algebras (see e.g. in \cite{AS}, \cite{ATVB}, \cite{RB4}, \cite{SmSt}) but, thanks to the analysis of \cite{RB3}, they admit the simplified formulation described here in the case of homogeneous algebras. For the latter ones we refer to \cite{BD-VW} the notations of which are used throughout. In the following $N$ is an integer with $N\geq 2$ and $\mathbb K$ is the ground field which is assumed to be algebraically closed and of characteristic zero.\\

A {\sl homogeneous algebra of degree $N$} or {\sl $N$-homogeneous algebra} is an algebra of the form 
\[
\cala = A(E,R)=T(E)/(R)
\]
where $E$ is a finite-dimensional vector space, $T(E)$ is the tensor algebra of $E$ and $(R)$ is the two-sided ideal of $T(E)$ generated by a vector subspace $R$ of $E^{\otimes^N}$. Such an algebra $\cala$ is naturally a $\mathbb N$-graded algebra $\cala=\oplus_n \cala_n$ with $\cala_0=\mathbb K$ and $\cala_1=E$, i.e. it is connected and generated in degree one.\\

Given a $N$-homogeneous algebra $\cala=A(E,R)$ one defines {\sl its dual} $\cala^!$ to be \cite{BD-VW} the $N$-homogeneous algebra $\cala^!=A(E^\ast,R^\perp)$ where $E^\ast$ is the dual vector space of $E$ and where $R^\perp \subset E^{\ast\otimes^N}=(E^{\otimes^N})^\ast$ is the annihilator of $R$.\\

To a $N$-homogeneous algebra $\cala=A(E,R)$ are canonically associated two dual $N$-complexes : The chain $N$-complex of left $\cala$-modules $K(\cala)$ and the cochain $N$-complex of right $\cala$-modules  $L(\cala)$, \cite{BD-VW}. One has $K(\cala)=\oplus_n K_n(\cala)$ with $K_n(\cala)=\cala\otimes \cala^{!\ast}_n$ and $L(\cala)=\oplus_n L^n(\cala)$ with $L^n(\cala)=\cala^!_n\otimes \cala$ where $\cala^{!\ast}_n\subset E^{\otimes^n}$ is the dual vector space of $\cala^!_n$ which is by definition a quotient of $E^{\ast\otimes^n}=(E^{\otimes^n})^\ast$. The $N$-differential $d:K_{n+1}(\cala)\rightarrow K_n(\cala)$ of $K(\cala)$ is induced by the left $\cala$-module homomorphism of $\cala\otimes E^{\otimes^{n+1}}$ into $\cala\otimes E^{\otimes^n}$ defined by
$a\otimes (e_0\otimes e_1\otimes \dots \otimes e_n)\mapsto (ae_0)\otimes (e_1\otimes \dots \otimes e_n)$
while the $N$-differential of $L(\cala)=\cala^!\otimes \cala$ is the left multiplication by the element $\xi_I$ of $E^\ast \otimes E\subset \cala^!\otimes \cala$ corresponding to the identity mapping $I$ of $E$ onto itself, i.e. $\xi_I =\theta^i\otimes e_i$ for any basis $(e_i)$ of $E$ with dual basis $(\theta^i)$. It is clear that these are $N$-differentials, i.e. satisfy $d^N=0$, and are (left or right) $\cala$-linear mappings. With obvious identifications, one passes from $K(\cala)$ to $L(\cala)$ by applying the functor $\Hom_\cala(\bullet,\cala)$ to each $\cala$-module $K_n(\cala)$. \\

Given a $N$-complex one can obtain from it ordinary complexes (2-complexes) called contractions by putting together alternatively $p$ or $N-p$ arrows $d$ of the $N$-complex (see e.g. \cite{D-V2}, \cite{BD-VW}). If $K$ is a chain $N$-complex these are the chain complexes $C_{m,p}(K)$
\[
\dots \stackrel{d^{N-m}}{\rightarrow} K_{Nr+p} \stackrel{d^m}{\rightarrow} K_{Nr+p-m} \stackrel{d^{N-m}}{\rightarrow} \dots
\]
and if $L$ is a cochain $N$-complex these are the cochain complexes $C_{m,p}(L)$
\[
\dots \stackrel{d^{N-m}}{\rightarrow} L^{Nr+p} \stackrel{d^m}{\rightarrow} L^{Nr+p+m} \stackrel{d^{N-m}}{\rightarrow} \dots
\]
where all the possibilities are covered by the conditions $0\leq p\leq N-2$ and $p+1\leq m\leq N-1$. In the case $N=2$ one has a complex which is of course its unique contraction $(C_{1,0})$. This is the case for $K(\cala)$ and $L(\cala)$ when $\cala$ is a quadratic algebra that is for $N=2$. In this latter case $K(\cala)$ is the {\sl Koszul complex} the acyclicity in positive degrees of which characterizes the {\sl quadratic Koszul algebras},  a class of very regular algebras which contains the algebras of polynomials.\\

In the case of a $N$-homogeneous algebra $\cala$ with $N\geq 3$, the right generalization of the above Koszul complex is the complex $C_{N-1,0}=C_{N-1,0}(K(\cala))$. Indeed as shown in \cite{BD-VW} it is the only contraction of $K(\cala)$ the acyclicity of which in positive degrees does not lead to a trivial class of algebras and furthermore it coincides with the Koszul complex of \cite{RB3} which was shown there to be a good generalization for $N\geq 3$ of the usual one (i.e. $K(\cala))$ in the quadratic case. Accordingly this complex $C_{N-1,0}(K(\cala))$ of left $\cala$-modules is refered to as the {\sl Koszul complex} of $\cala$ and $\cala$ is said to be a {\sl Koszul algebra} whenever this complex is acyclic in positive degrees. Since in degree 0 the homology of the Koszul complex of $\cala$ is the trivial left $\cala$-module $\mathbb K$, if $\cala$ is a Koszul algebra the Koszul complex of $\cala$ is a resolution of the trivial left $\cala$-module $\mathbb K$.\\

Given a positively graded complex $\calc=\oplus_{n\geq 0} \calc_n$ of chains or of cochains, $\calc$ is said to have {\sl finite extension} $D\in \mathbb N$ if $\calc_n=0$ for $n >D$ and $\calc_D\not=0$. The dual of the Koszul complex $C_{N-1,0}$ of a $N$-homogeneous algebra $\cala$ is the complex of right $\cala$-modules $C'_{N-1,0}=C_{1,0}(L(\cala))$. By dual of the chain complex of left $\cala$-modules $C_{N-1,0}$ we here mean the cochain complex $C'_{N-1,0}$ of right $\cala$-modules obtained by applying the functor $\Hom_\cala(\bullet,\cala)$ to each left $\cala$-module  of $C_{N-1,0}$. If the Koszul complex of $\cala$ has finite extension $D$, the same is true for its dual $C'_{N-1,0}$. If $\cala$ is a Koszul algebra with Koszul complex $C_{N-1,0}$ of finite extension $D$ then $\cala$ is of {\sl global dimension} $D$. This implies in particular \cite{RB3} that $\cala$ has Hochschild dimension $D$. If furthermore the dual complex $C'_{N-1,0}$ gives a resolution of the trivial right $\cala$-module $\mathbb K$, i.e. if $H^k(C'_{N-1,0})=0$ for $k<D$ and $H^D(C'_{N-1,0})\simeq \mathbb K$, then $\cala$ is said to be Gorenstein. This is clearly a particular case of a sort of (generalization of) Poincar\'e duality property.

\end{document}